\documentclass[11pt]{article}
\usepackage{latexsym}
\usepackage{amssymb}
\newlength{\mytopmargin}
\newlength{\myleftmargin}
\newtheorem{lemma}{Lemma}[section]

\newtheorem{thm}[lemma]{Theorem}

\newtheorem{prop}[lemma]{Proposition}

\newcommand{\zz}{\mathbb{Z}}

\newcommand{\inn}[2]{\left\langle #1, #2 \right\rangle}
\newcommand{\harg}[3]{\left[ {#1 \atop #2} ; #3 \right]}

\begin{document}
\noindent
\begin{center}{ \Large\bf Transformation formulae for 
multivariable \\ basic hypergeometric series}
\end{center}
\vspace{5mm}
\noindent
\begin{center} T.~H.~Baker and P.~J.~Forrester \\[3mm]
{\it Department of Mathematics, University of Melbourne, \\
Parkville, Victoria 3052, Australia} 
\end{center}
\vspace{.5cm}

\begin{quote} We study multivariable (bilateral) basic hypergeometric 
series associated with (type $A$) Macdonald polynomials. We derive
several transformation and summation properties for such series including
analogues of Heine's ${}_2\phi_1$ transformation, the $q$-Pfaff-Kummer
and Euler transformations, the $q$-Saalsch\"utz summation formula and
Sear's transformation for terminating, balanced ${}_4\phi_3$ series.
For bilateral series, we rederive Kaneko's analogue of the ${}_1\psi_1$
summation formula and give multivariable extensions of Bailey's
${}_2\psi_2$ transformations.
\end{quote}
%\noindent {\bf Mathematics Subject Classification:}\quad 33D80
\begin{center}
{\it Dedicated to Dick Askey on the occasion of his 65'th birthday.}
\end{center}
 
\setcounter{section}{-1}
\section{Introduction}
Multivariable basic hypergeometric series of the type studied in this 
paper were first introduced by Kaneko and Macdonald \cite{kaneko96a,macunp1}.
They are defined as
\begin{equation} \label{basic.def}
{}_r\Phi_s\harg{a_1,\ldots,a_r}{b_1,\ldots,b_s}{z} :=
\sum_{\lambda} \left((-1)^{|\lambda|}q^{n(\lambda')}\right)^{s+1-r}
\frac{(a_1)_{\lambda}\cdots(a_r)_{\lambda}}
{(b_1)_{\lambda}\cdots(b_s)_{\lambda}h'_{\lambda}} \; P_{\lambda}(z;q,t)
\end{equation}
where all quantities are defined in Section \ref{two}.

The $q=1$ case had been studied previously by Yan \cite{yan92b}
where several important properties of the hypergeometric series
${}_p {\cal F}_q$ were studied including the Gauss 
(${}_2 {\cal F}_1$) and Kummer
(${}_1 {\cal F}_1$) summation formulae, the Pfaff-Kummer and Euler
transformations and integral representations were derived (see also
\cite{opdam93a} for the Gauss formula with respect to arbitrary root
systems).

Kaneko \cite{kaneko96a} considered a generalized $q$-Selberg integral
dependent on parameters $x_1,\ldots,x_m$ and derived a set of 
$m$ $q$-difference equations satisfied by such an integral. He then
showed that the basic hypergeometric series ${}_2\Phi_1$ defined by
(\ref{basic.def}) was the unique solution (satisfying
certain properties) of such a system. He was thus able to derive an
integral representation for this ${}_2\Phi_1$ series and hence give 
an alternative proof of the $q$-Selberg integral \cite{askey80a}
(see also \cite{kad88a,zeil89a,kad94a}). He additionally derived the
$q$-analogue of the Gauss formula and another integral formula,
the constant term version of which was presented in 
\cite[Theorem 4]{kaneko97b}.
In a subsequent work \cite{kaneko97a}, Kaneko introduced a multivariable
analogue of the bilateral basic hypergeometric series ${}_r\Phi_s$
and derived an analogue of Ramanujan's ${}_1\psi_1$ summation formula,
along with a multivariable version of the Jacobi triple product 
identity (which is a limiting case - see also \cite{kaneko96b}).

Independently, Macdonald in his unpublished notes \cite{macunp1} carried
out (among other things) a similar program to Yan in the $q=1$ case along 
with a multivariable version of the Saalsch\"utz summation formula (the
summation of a balanced ${}_3 {\cal F}_2$ with unit argument), while for 
general $q$ he derived the integral representation for the ${}_2\Phi_1$ 
series, and $q$-analogues of the Gauss and Saalsch\"utz formulae.

The aim of the present work is to supplement some of the existing
knowledge with some new transformation and summation formulae in this
multivariable setting (including the $q$-analogue of the Pfaff-Kummer
transformation formulae, Heine's transformation formula for the 
${}_2\Phi_1$ series, Sear's transformation for
terminating, balanced ${}_4\Phi_3$ series, and various summation and
transformation formulae for bilateral ${}_2\Psi_2$ series),
as well as providing alternative derivations of known results. For a
comprehensive review of summation and transformation formulae for
basic hypergeometric series in the one-variable case see Gasper and
Rahmen's book \cite{gasperbk}.

We note here that many of the summation and transformation formulae 
involving multivariable basic hypergeometric series are only valid when
the argument is specialized to $zt^{\delta}:=(z,zt,zt^2,\ldots,zt^{n-1})$ 
(the exceptions being the $q$-binomial theorem, the Euler transformation 
of the ${}_2\Phi_1$ series and the ${}_1\Psi_1$ summation formula).

Finally, we point out that many of the formulae presented here can be 
derived in the case $q=t$ (the Schur polynomial case) as special cases
of the very general formulae of Milne and co-workers (see 
\cite{miln97a} and referenced therein), or from those found in the
works of Krattenthaler, Gustafson and Schlosser \cite{krat97a,schloss96a}.
Also, summation formulae for hypergeometric systems associated with the
$BC$ root system have been considered by Van Diejen \cite{vandiej97a}
(see also \cite{vandiej97b}).

The plan of the paper is as follows. In Section 1 we set out the basic 
facts about Macdonald polynomials we shall require. Section 2 exhibits
a multivariable extension of Heine's ${}_2\phi_1$ transformations. The
Euler transformation is derived in Section 3 from the defining
difference equations for the ${}_2\Phi_1$ series. Section 4 treats the
$q$-Saalsch\"utz formula while the $q$-Pfaff-Kummer transformation
for ${}_2\Phi_1$ and Sear's transformation for terminating, balanced
${}_4\Phi_3$ series are discussed in Section 5. Bilateral series are
studied in Section 6 where the ${}_1\Psi_1$ summation formula is 
derived using a shifted version of the multivariable Gauss summation
formula, and various ${}_2\Psi_2$ transformations are derived using
two different methods. 

\setcounter{equation}{0}
\section{Notations} \label{two}

The Macdonald polynomials $P_{\lambda}(x;q,t)$, $x:=(x_1,\ldots,x_n)$ 
(often abbreviated
to $P_{\lambda}(x)$ or just $P_{\lambda}$ when the context is clear) are
defined as the unique symmetric polynomials having the expansion
$$
P_{\lambda}(x) = m_{\lambda}(x) + \sum_{\mu<\lambda} c_{\lambda\mu}
m_{\mu}(x)
$$
(where $<$ denotes the dominance order on partitions and $m_{\lambda}(x)$
denotes the monomial symmetric function), which form
an orthogonal basis of symmetric functions with respect to the inner
product
\begin{equation} \label{inner}
\inn{f}{g} := \frac{1}{n!}\;{\rm C.T.} \left( f(x)\, g(x^{-1}) 
\Delta_q(x) \right)
\end{equation}
where
\begin{equation} \label{dedo}
\Delta_q(x) := \prod_{1\leq i<j \leq n} (x_i/x_j;q)_k (x_j/x_i;q)_k
\end{equation}
Here ${\rm C.T.}$ stand for ``the constant term of'' in the case where
$t=q^k$ (and hence the quantity in brackets in (\ref{inner}) is a
Laurent polynomial) or the corresponding trigonometric integral otherwise,
and
$$
(x;q)_a := \frac{(x)_{\infty}}{(xq^a;q)_{\infty}}, \hspace{3cm}
(u;q)_{\infty} := \prod_{i=0}^{\infty} (1-uq^i) .
$$
Define the following quantities
\begin{equation} \label{libro}
h_{\lambda} = \prod_{s\in\lambda} (1-q^{a(s)}t^{l(s)+1}) \qquad
h'_{\lambda} = \prod_{s\in\lambda} (1-q^{a(s)+1}t^{l(s)}) 
\end{equation}
where $a(s)$ (respectively $l(s)$) denotes the arm-length (resp. leg-length)
of the node $s$ in the diagram of $\lambda$,
along with the generalized shifted $q$-factorial
$(a)_{\lambda}:=t^{n(\lambda)}\,\prod_{i=1}^n (at^{1-i};q)_{\lambda_i}$.
Here $n(\lambda) = \sum_i (i-1)\lambda_i = \sum_i \lambda'_i(\lambda'_i-1
)/2$.
Two fundamental results of Macdonald polynomials we shall use, both of
which can be found in \cite{mac}, is their specialization formula
\begin{equation} \label{eval}
P_{\lambda}(1,t,\ldots,t^{n-1}) = \frac{(t^n)_{\lambda}}{h_{\lambda}}
\end{equation}
and their normalization with respect to the inner product (\ref{inner})
\begin{equation} \label{norm}
\inn{P_{\lambda}}{P_{\lambda}} = \frac{h'_{\lambda}}{h_{\lambda}}
\frac{(t^n)_{\lambda}}{(qt^{n-1})_{\lambda}}\: \inn{1}{1}
\end{equation}
where
\begin{equation} \label{ground}
\inn{1}{1} = \frac{1}{n!}\;{\rm C.T.}\left( \prod_{1\leq i < j \leq n}
(x_i/x_j;q)_k (x_j/x_i;q)_k \right) = \prod_{i=1}^n \frac{(q;q)_{ik-1}}
{(q;q)_{k-1} (q;q)_{(i-1)k}}
\end{equation}
 
Further, following Macdonald \cite{mac}, for a 
partition $\lambda$, let $u_{\lambda}$ 
denote the evaluation map on polynomials in $n$ variables which sets 
$x_i=q^{\lambda_i}t^{n-i}$. There are no nice expressions for $u_{\lambda}
(P_{\mu})$ except in the cases $\lambda=0$ (since $u_0(P_{\mu})=
P_{\mu}(1,t,\ldots,t^{n-1})$ given above in (\ref{eval})) or more
generally $\lambda=(m^n)$ (since $u_{(m^n)}(P_{\mu}) = q^m u_0(P_{\mu})$).
Nonetheless, there is a useful symmetry property for $u_{\lambda}
(P_{\mu})$ in general, which reads
\begin{equation} \label{symm.mac}
\frac{u_{\lambda}(P_{\mu})}{u_0(P_{\mu})} =
\frac{u_{\mu}(P_{\lambda})}{u_0(P_{\lambda})} 
\end{equation}

One of the fundamental results in the theory of basic hypergeometri series
is the summation formula for a ${}_1\phi_0$ series, also known as
the $q$-binomial theorem
$$
\sum_{n\geq 0} \frac{(a;q)_n}{(q;q)_n} \; z^n =
\frac{(az;q)_{\infty}}{(z;q)_{\infty}} 
$$
This has the following important generalization (see e.g. 
\cite[Thm 3.5]{kaneko96a})
\begin{equation} \label{q.bin}
{}_1\Phi_0\harg{a}{-}{z} = \prod_{i=1}^n \frac{(a z_i;q)_{\infty}}
{(z_i;q)_{\infty}}
\end{equation}
which will often be of use below.

\setcounter{equation}{0}
\section{Heine's transformation formula}
A fundamental transformation formula for ${}_2\phi_1$ series is
Heine's transformation formula
\begin{equation} \label{heine.1}
{}_2\phi_1\harg{a,b}{c}{z} = \frac{(b)_{\infty} (az)_{\infty}}
{(c)_{\infty} (z)_{\infty}} \: {}_2\phi_1 \harg{c/b,z}{az}{b}
\end{equation}
The utility of this formula is two-fold: iterate it twice to get the
Euler transformation formula (see (\ref{euler.1}); set 
$z=c/ab$ and the ${}_2\phi_1$
on the right hand side of (\ref{heine.1}) becomes a ${}_1\phi_0$ which
can be summed by the $q$-binomial theorem resulting in Gauss's summation
formula 
$$
{}_2\phi_1\harg{a,b}{c}{\frac{c}{ab}} = \frac{(c/a)_{\infty}
(c/b)_{\infty}}{(c)_{\infty} (c/ab)_{\infty}}
$$
The multivariable version of Heine's formula appears
as follows
\begin{equation} \label{heine.n}
{}_2\Phi_1\harg{a,b}{c}{z t^{\delta}}=\prod_{i=1}^n \frac{(bt^{1-i})_{\infty}
(azt^{n-i})_{\infty}}{(ct^{1-i})_{\infty} (zt^{n-i})_{\infty}}\:
{}_2\Phi_1\harg{c/b,zt^{n-1}}{azt^{n-1}}{bt^{1-n} t^{\delta}}
\end{equation}
(as mentioned in the introduction, the notation $z t^{\delta}$ stands for
the argument $(z,zt,\ldots,zt^{n-1})$).
The proof of this identity follows in the same manner as the one-variable
case, with the aid of the symmetry property (\ref{symm.mac}).
To see this, first note that setting $z_i=b t^{1-i}q^{\lambda_i}$ and
$a=c/b$ in the $q$-binomial theorem (\ref{q.bin}) yields
\begin{equation} \label{dokusho}
\prod_i \frac{(c t^{1-i}q^{\lambda_i})_{\infty}}
{(b t^{1-i}q^{\lambda_i})_{\infty}} = \sum_{\mu} \frac{(c/b)_{\mu}}
{h'_{\mu}}\: (b t^{1-n})^{|\mu|}\; u_{\lambda}(P_{\mu})
\end{equation}
Thus
\begin{eqnarray*}
{}_2\Phi_1\harg{a,b}{c}{z t^{\delta}} &=& 
\prod_{i=1}^n \frac{(b t^{1-i})_{\infty}}
{(c t^{1-i})_{\infty}}\sum_{\lambda} \frac{(a)_{\lambda}}{h'_{\lambda}}
P_{\lambda}(t^{\delta}) z^{|\lambda|} \sum_{\mu} 
\frac{(c/b)_{\mu}}
{h'_{\mu}}\: (b t^{1-n})^{|\mu|}\; u_{\lambda}(P_{\mu}) \\
&=&  \prod_{i=1}^n \frac{(b t^{1-i})_{\infty}} {(c t^{1-i})_{\infty}}
\sum_{\lambda,\mu} \frac{(a)_{\lambda} (c/b)_{\mu}}{h'_{\lambda}
h'_{\mu}} z^{|\lambda|} (bt^{1-n})^{|\mu|} u_0(P_{\mu}) u_{\mu}
(P_{\lambda}) \\
&=& \prod_{i=1}^n \frac{(b t^{1-i})_{\infty}(azt^{n-i})_{\infty}} 
{(c t^{1-i})_{\infty} (zt^{n-i})_{\infty}} \sum_{\mu}
\frac{(c/b)_{\mu} (z t^{n-1})_{\mu}}{(azt^{n-1})_{\mu}h'_{\mu}}
(bt^{1-n})^{|\mu|} u_0(P_{\mu})
\end{eqnarray*}
which gives the result. Here we have used (\ref{dokusho}), (\ref{symm.mac})
and (\ref{dokusho}) respectively.

As mentioned above, upon setting $z=c/(abt^{n-1})$ the r.h.s.
of (\ref{heine.n}) reduces to a ${}_1\Phi_0$ which can be summed through
the $q$-binomial theorem giving Gauss's formula
\begin{equation} \label{gauss.n}
{}_2\Phi_1 \harg{a,b}{c}{\frac{c}{abt^{n-1}} t^{\delta}} = \prod_{i=1}^n
\frac{(\frac{c}{b}t^{1-i})_{\infty}(\frac{c}{a}t^{1-i})_{\infty}}
{(\frac{c}{ab}t^{1-i})_{\infty}(ct^{1-i})_{\infty}}
\end{equation}
Also, we can iterate (\ref{heine.n}) twice and obtain a version of
Euler's transformation valid for the variables $zt^{\delta}$. In fact,
we shall see in the next section that there exists a multivariable
version of the Euler transformation, which is true 
for general argument $z:=(z_1,\ldots,z_n)$.

\setcounter{equation}{0}
\section{Euler transformation}
In the theory of the one-variable $q$-hypergeometric function ${}_2\phi_1$,
the Euler transformation reads
\begin{equation}\label{euler.1}
{}_2\phi_1\Big [ {a,\,b \atop c};x \Big ] =
{(abx/c)_\infty \over (x)_\infty} {}_2 \phi_1 
\Big [ {c/a,\,c/b \atop c};abx/c \Big ].
\end{equation}
This transformation generalizes naturally to the $n$-variable case.

\begin{prop}
With ${}_2\Phi_1$ defined by (\ref{basic.def}), we have
\begin{equation}\label{euler.2}
{}_2\Phi_1\Big [ {a,\,b \atop c};z \Big ] =
\prod_{i=1}^n {(abz_i/c)_\infty \over (z_i)_\infty} {}_2 \Phi_1 
\Big [ {c/a,\,c/b \atop c};abz/c \Big ].
\end{equation}
\end{prop}

The  $n$-variable Euler transformation (\ref{euler.2}) can be proved by
using a theorem of Kaneko \cite[Thm.~4.12]{kaneko96a} which characterizes
${}_2\Phi_1$ as the unique solution of a system of $q$-difference equations.
To state these equations, let $\tau_i$ denote the $q$-shift operator for
the variable $z_i$ so that $\tau_i f(z_1,\dots,z_n) = f(z_1,\dots,qz_i, \dots,
z_n)$, write
$$
A_i(z;t) := \prod_{l=1 \atop l \ne i}^n {tz_i - z_l \over z_i - z_l},
$$
and define the $q$-derivative by
\begin{equation}\label{qder}
{\partial \phi \over \partial_q z_i} = {(1-\tau_i) \phi \over
(1-q) z_i}.
\end{equation}
\begin{thm} \label{thm.k}
{\bf (Kaneko)} The multivariable basic $q$-hypergeometric series
${}_2\Phi_1\Big [ {a,\,b \atop c};z \Big ]$ as defined by (\ref{basic.def})
is the unique solution of the system of $q$-difference equations\footnote{This
is eqn.~(2.26) of \cite{kaneko96a} with the factor $(1-q)$ in the second 
term corrected to read $(1-t)$}
\begin{eqnarray}\label{qdif}
z_i(c-abqz_i)\tau_i(A_i(z;t)) {\partial^2 S \over \partial_q z_i^2}
+(1-t) \sum_{j=1 \atop j \ne i}^n {z_i z_j (c - abz_j) \over
qz_i - tz_j} \tau_i(A_j(z;t)) {\partial^2 S \over \partial_q z_i  
\partial_q z_j} \nonumber \\
+ \bigg \{ {t^{n-1} - c \over 1 - q} + {1 \over 1 -q} \Big (
(1-a)(1-b) t^{n-1} - (t^{n-1} - abq) \Big ) z_i \bigg \}
{\partial S \over \partial_qz_i} \nonumber \\
+ {1 - t \over 1 - q} \bigg \{ {1 - \tau_i(A_i(z;t)) \over 1 - t}
(c - abqz_i) {\partial S \over \partial_q z_i} 
-  \sum_{j=1 \atop j \ne i}^n { z_j (c - abz_j) \over
qz_i - tz_j} \tau_i(A_j(z;t)) {\partial S \over   
\partial_q z_j} \bigg \} \nonumber \\
- {(1-a)(1-b) t^{n-1} \over (1-q)^2} S = 0 \qquad (i=1,\dots,n),
\end{eqnarray}
subject to the conditions that $S(z)$ is a symmetric function in
$z_1,\dots,z_n$ and $S(z)$ is analytic at the origin with $S(0)=1$.
\end{thm}

In light of Theorem \ref{thm.k}, the Euler transformation (\ref{euler.2})
can be proved by verifying that the substitution
$$
S = \prod_{i=1}^n {(abz_i/c)_\infty \over (z_i)_\infty} \; U
$$
into the $q$-difference equation (\ref{qdif}) implies that $U$ satisfies
the same $q$-difference equation (\ref{qdif}) with the replacements
\begin{equation}\label{replace}
a \mapsto {c \over a}, \quad b \mapsto {c \over b}, \quad
c \mapsto c, \quad z_i \mapsto {abz_i \over c}.
\end{equation}
To perform this task we will require the $q$-differentiation product
rule
\begin{equation}\label{pr}
{ \partial(\phi \psi) \over \partial_q z_i} =
{ \partial\phi  \over \partial_q z_i} \psi +
(\tau_i \phi) { \partial\psi \over \partial_q z_i},
\end{equation}
as well as the special summation formula \cite[eqn.~(2.11)]{kaneko96a}
\begin{equation}\label{sum}
\sum_{j=1 \atop j \ne i}^n
{z_j A_j(z;t) \over tz_j - z_i} = {A_i(z;t) - t^{n-1} \over 1 - t}.
\end{equation}

\vspace{.2cm}
\noindent
{\bf Proof of Proposition 3.1} \\
Let us now give some details of the required calculation. With
$$
P_i := {(abz_i/c)_\infty \over (z_i)_\infty}, \qquad
P_i'' := {(abq^2z_i/c)_\infty \over (z_i)_\infty},
$$
straightforward use of (\ref{qder}) and (\ref{pr}) shows that
\begin{eqnarray}
{\partial \over \partial_q z_i}(P_i U) & = &
{1 - ab/c \over 1 - q} (1 - abqz_i/c) P_i'' U +
(1 - z_i)(1 - abqz_i/c) P_i'' {\partial U \over \partial_q z_i} 
\label{1.one} \\{\partial^2 \over \partial_q z_i^2}(P_i U) & = &
{(1 - ab/c)(1-abq/c)  \over (1 - q)^2}  P_i'' U +
{(1 - ab/c)(1-z_i)  \over 1 - q}P_i'' {\partial U \over \partial_q z_i}
\nonumber  \\
&& +{q \over 1 - q}(1 - ab/c)(1 - z_i) P_i'' {\partial U \over \partial_q
z_i} + (1 - z_i) (1 - qz_i) P_i'' {\partial^2 U \over \partial_q z_i^2}
\end{eqnarray}
Double use of (\ref{1.one}) allows us to write down a similar formula for
$$
{\partial^2  \over \partial_q z_i \partial_q z_j}(P_iP_jU).
$$

Next we use these formulas to rewrite the derivatives in
(\ref{qdif}) with $S = (\prod_{i=1}^n P_i) \, U$, and proceed to collect
together terms so that the essential structure of (\ref{qdif}) is
maintained.
Thus, for example, we must collect together all terms proportional to $U$,
and show that apart from a common factor (which is eventually canceled out)
the coefficient of $S$ in (\ref{qdif}), with the substitutions (\ref{replace}),
results. These terms are
\begin{eqnarray}\label{U}
{(1 - abqz_i/c) \over (1 - q)^2} \Big ( \prod_{l=1 \atop l \ne i}^n
P_l \Big ) P_i'' \bigg [
cz_i \tau_i(A_i(z;t)) (1 - ab/c)(1-abq/c) \nonumber \\
+ (1 - t) c (1 - ab/c)^2 z_i \sum_{j=1 \atop j \ne i}^n
{z_j \over qz_i - tz_j} \tau_i(A_j(z;t)) \nonumber \\
+ (1-ab/c) \Big \{ t^{n-1} - c + \Big (
(1-a)(1-b)t^{n-1} - (t^{n-1} - abq) \Big )z_i \Big \} \nonumber \\
+(1-t)(1-ab/c)\Big \{
{1 - \tau_i(A_i(z;t)) \over 1 - t}(c - abqz_i) - c  \sum_{j=1 \atop j \ne i}^n
{z_j \over qz_i - tz_j} \tau_i(A_j(z;t)) \Big \} \nonumber \\
- (1-a)(1-b) t^{n-1} (1 - abz_i/c) \bigg ].
\end{eqnarray}
Noting that
$$
 \sum_{j=1 \atop j \ne i}^n
{z_j \over qz_i - tz_j} \tau_i(A_j(z;t)) =
\tau_i \bigg (  \sum_{j=1 \atop j \ne i}^n
{z_j \over z_i - tz_j} A_j(z;t) \bigg ),
$$
we see that the summation can be evaluated according to (\ref{sum}).
The expression (\ref{U}) then greatly simplifies, reducing down to
\begin{equation}\label{F}
\bigg ((1-abqz_i/c)(1-z_i) \Big ( \prod_{l=1 \atop l \ne i}^n
P_l \Big ) P_i'' \bigg ) \bigg \{
-{(1-a/c)(1-b/c) t^{n-1} \over (1-q)^2} \bigg \}.
\end{equation}
The factor in the curly brackets is precisely what is required to obtain
(\ref{qdif}) with the replacements (\ref{replace}).

Similar calculations to simplify the terms proportional to the
$q$-derivatives show that the precise $q$-difference equation
(\ref{qdif}) results with $S=U$ and the replacements (\ref{replace});
the term in the parenthesis in (\ref{F}) is a common factor which
cancels out. \hfill$\Box$

\setcounter{equation}{0}
\section{$q$-Saalsch\"utz summation formula}

A summation formula for basic hypergeometric series used in many 
contexts is the sum of a terminating balanced ${}_3\phi_2$ series 
($q$-Saalsch\"utz formula)
$$
{}_3\phi_2\harg{q^{-N},b,c}{e,q^{1-N}bc/e}{q} = \frac{(e/b)_N (e/c)_N}
{(e)_N (e/bc)_N}
$$
This formula follows directly by comparing coefficients of $z^N$ on
both sides of the Euler transformation formula (\ref{euler.1}) 
(the products being expanded with the $q$-binomial theorem).
This can be generalized to a multi-variable setting following the proof
of Macdonald \cite{macunp1} in the $q=1$ case. First, we need some 
preliminary results. 

Given a partition $\lambda$ of length $\ell(\lambda)\leq n$
let $\hat{\lambda}$ to be the complement of $\lambda$ in the rectangle
$(N^n)$. That is, $\hat{\lambda}_i = N - \lambda_{n+1-i}$. Then it follows
from the various definitions that
\begin{equation} \label{ototo}
(a)_{\hat{\lambda}} = t^{n(\lambda)} q^{n(\lambda')-(N-1)|\lambda|}
\frac{(a)_{(N^n)}}{(-a)^{|\lambda|} (t^{n-1}q^{1-N}/a)_{\lambda}}
\end{equation}

Next, let $f^{\lambda}_{\mu\nu}(q,t)$ denote the coefficient of $P_{\lambda}$
in the expansion of $P_{\mu}\,P_{\nu}$. The following result is due to
Macdonald \cite{macunp1} but we include the proof for completeness,
\begin{lemma}
We have
\begin{equation} \label{imoto}
f^{(N^n)}_{\mu\nu}(q,t) = \frac{\inn{P_{\mu}}{P_{\mu}}}
{\inn{1}{1}} \; \delta_{\hat{\mu} \nu} = \frac{(t^n)_{\mu} h'_{\mu}}
{(qt^{n-1})_{\mu} h_{\mu}}\;\delta_{\hat{\mu} \nu}
\end{equation}
\end{lemma}
{\it Proof.}\quad Writing $|x|:= x_1\cdots x_n$, first note that
\begin{equation} \label{sombrero}
P_{\hat{\lambda}}(x) = |x|^N\: P_{\lambda}(x^{-1})
\end{equation}
which follows from the characterization of the Macdonald polynomials
explained in Section \ref{two} (clearly both sides of the above 
equation have leading order term $m_{\hat{\lambda}}(x)$ and also \linebreak
$\inn{|x|^N P_{\lambda}(x^{-1})}{|x|^N P_{\mu}(x^{-1})} = 
\inn{P_{\mu}}{P_{\lambda}}$, yielding the requisite orthogonality
 -- see also \cite{vandiej97b} on this point).
It thus follows from (\ref{sombrero}) that
\begin{eqnarray*}
\inn{P_{\lambda}(x)}{P_{\mu}(x)} = \inn{P_{\lambda}(x)\,P_{\mu}(x^{-1})}{1}
= \inn{P_{\lambda}(x) P_{\hat{\mu}}(x)}{|x|^N} \\
=\sum_{\sigma} f^{\sigma}_{\lambda\hat{\mu}}(q,t) \inn{P_{\sigma}(x)}
{P_{(N^n)}(x)} = f^{(N^n)}_{\lambda\hat{\mu}}(q,t) \inn{1}{1}
\end{eqnarray*}
The result now follows from the orthogonality of Macdonald polynomials
and the expression (\ref{norm}) for their normalization.\hfill $\Box$

The last ingredient needed is the following expression for $h'_{\lambda}$
due to Kaneko \cite{kaneko96b}
\begin{equation} \label{lapiz}
h'_{\lambda} = (q)_{\infty}^n \prod_i \frac{1}{(q^{\lambda_i+1}t^{n-i})_{
\infty}} \prod_{i<j} \frac{(q^{\lambda_i-\lambda_j+1}t^{j-i})_{\infty}}
{(q^{\lambda_i-\lambda_j+1}t^{j-i-1})_{\infty}}
\end{equation}
from which one obtains a relation we shall need presently
\begin{equation} \label{shumi}
\frac{h'_{\hat{\lambda}}}{h'_{\lambda}} = (-1)^{|\lambda|} 
t^{n(\lambda)-n((N^n))} q^{n(\lambda')-N|\lambda|} \frac{
(qt^{n-1})_{(N^n)}}{(q^{-N})_{\lambda} (qt^{n-1})_{\lambda}}
\end{equation}

Beginning with the Euler transformation formula (\ref{euler.2}), expand 
the products occurring there with the $q$-binomial theorem and compare 
coefficients of $P_{\lambda}(z)$ on both sides to obtain
\begin{equation} \label{saal.n}
\frac{(a)_{\lambda} (b)_{\lambda}}{(c)_{\lambda} h'_{\lambda}} =
\sum_{\mu, \nu} \frac{(c/a)_{\mu} (c/b)_{\mu} (ab/c)_{\nu}}
{(c)_{\mu} h'_{\mu} h'_{\nu}} \left(\frac{ab}{c}\right)^{|\mu|}
f^{\lambda}_{\mu\nu}(q,t)
\end{equation}
The point here is that in the case when $\lambda=(N^n)$, the sum on
the right hand side can be expressed as a balanced ${}_3\Phi_2$ series
(the series ${}_r\Phi_s$ in (\ref{basic.def}) is said to be {\it balanced}
if $a_1\cdots a_r qt^{n-1} = b_1\cdots b_s$). Indeed after using the 
expression (\ref{imoto}) for $f^{(N^n)}_{\mu\nu}$ and expressing
all quantities labelled by $\hat{\mu}$ in terms of those of $\mu$
by means of (\ref{ototo}) and (\ref{shumi}) the resulting sum over
$\mu$ can be simplified to yield 
\begin{equation} \label{saal.n2}
\frac{(a)_{(N^n)} (b)_{(N^n)}}{(c)_{(N^n)} (ab/c)_{(N^n)}} =
{}_3\Phi_2\harg{q^{-N},c/a,c/b}{c,\frac{c}{ab}t^{n-1}q^{1-N}}{q\,t^{\delta}}
\end{equation}
which is a multivariable version of the $q$-Saalsch\"utz formula.

Note that by taking the limit $a\rightarrow 0$ in (\ref{saal.n}) 
using the fact that $\lim_{x\rightarrow \infty} (x)_{\mu}/x^{|\mu|} =
(-1)^{|\mu|} q^{n(\mu')}$, we have (after shifting $b\rightarrow c/b$)
\begin{equation} \label{chu}
t^{n(\lambda)} \frac{(c/b)_{\lambda}}{(c)_{\lambda} h'_{\lambda}} =
\sum_{\mu,\nu} t^{n(\nu)} q^{n(\mu')} \frac{(b)_{\mu}}
{(c)_{\mu} h'_{\mu} h'_{\nu}} \left(-\frac{c}{b}\right)^{|\mu|}
f^{\lambda}_{\mu\nu}(q,t)
\end{equation}
which is a generalization of the Chu-Vandermonde formula (especially
in the case $\lambda = (N^n)$ when the right hand side can be expressed
as a terminating ${}_2\Phi_1$) we shall have subsequent need for.

\setcounter{equation}{0}
\section{Other transformation formulae}

One of the simplest transformations where the simple $n=1$ proof must be
rethought is the Pfaff-Kummer transformation
\begin{equation} \label{kumm.q}
{}_2\phi_1\harg{a,b}{c}{z} = \frac{(a z;q)_{\infty}}{(z;q)_{\infty}}\:
{}_2\phi_2\harg{a,c/b}{c,az}{bz}
\end{equation}
This is a $q$-analogue of the classical relation
\begin{equation} \label{kumm}
{}_2 F_1\harg{a,b}{c}{z} = (1-z)^{-a}\,{}_2 F_1\harg{a,c-b}{c}{\frac{z}
{z-1}}
\end{equation}
Indeed, it can be shown using the ${}_2 {\cal F}_1$ integral representation 
\cite{macunp1} or the defining set of $m$ partial differential equations
\cite{yan92b,kaneko93a} that the multivariable analogue of (\ref{kumm}) 
holds in its full generality:
$$
{}_2 {\cal F}_1\harg{a,b}{c}{z} = \prod_{i=1}^n (1-z_i)^{-a} 
{}_2 {\cal F}_1 \harg{a,c-b}{c}{\frac{z}{z-1}}
$$
This is not possible in the multivariable version of (\ref{kumm.q}); however
for the restricted argument $zt^{\delta}$, such a formula can be
derived via a suitable modification of the argument in \cite{gasperbk}.

The key step is to establish the identity
\begin{equation} \label{chichi}
\prod_{i=1}^n \frac{(azt^{n-i})_{\infty}}{(zt^{n-i})_{\infty}}
\frac{(a)_{\mu}}{(azt^{n-1})_{\mu}} u_0(P_{\mu}) =
\sum_{\lambda,\nu} t^{b(\nu)-b(\lambda)} \frac{(a)_{\lambda}}{h'_{\nu}}
u_0(P_{\lambda}) z^{|\nu|} f^{\lambda}_{\mu\nu} ,
\end{equation}
for then it follows from the Chu-Vandermonde formula (\ref{chu}) that
\begin{eqnarray}
{}_2\Phi_1\harg{a,b}{c}{z t^{\delta}} &=& \sum_{\lambda,\mu,\nu} 
(-b)^{|\mu|} q^{b(\mu')} t^{b(\nu) - b(\lambda)}
\frac{(c/b)_{\mu}}{(c)_{\mu} h'_{\mu}} \; \frac{(a)_{\lambda}}
{h'_{\lambda}} u_0(P_{\lambda}) z^{|\lambda|} f^{\lambda}_{\mu\nu}\;
z^{|\mu|+|\nu|} \nonumber\\
&=& \prod_{i=1}^n \frac{(a z t^{n-i};q)_{\infty}}
{(z t^{n-i};q)_{\infty}}\:
{}_2\Phi_2\harg{a,c/b}{c,azt^{n-1}}{bz t^{\delta}}
\end{eqnarray}
which is the multi-variable version of the Pfaff-Kummer transformation
formula. It thus remains to prove (\ref{chichi}). This identity is in
fact a special case of the more general identity
\begin{lemma} \label{gauss.shift}
$$
\sum_{\sigma,\mu} \frac{(c)_{\sigma} z^{|\sigma|}}{h'_{\sigma}}\;
\frac{(a)_{\mu} x^{|\mu|} P_{\mu}(t^{\delta})}{(b)_{\mu}}\;
f^{\mu}_{\lambda\sigma}(q,t) =
\prod_{i=1}^n \frac{(czxt^{n-i})_{\infty} 
(azxt^{n-i})_{\infty}}{ (bt^{1-i})_{\infty} (zxt^{n-i})_{\infty} }
\; \frac{(ax)_{\lambda}}{(azxt^{n-1})_{\lambda}}
x^{|\lambda|} P_{\lambda}(t^{\delta})
$$
provided $b=aczxt^{n-1}$. 
\end{lemma}
{\it Proof.}\quad Let $S$ denote the summation side of the above
equation. Note from the definition of the generalized $q$-factorials
and the $q$-binomial theorem we can write
$$
\frac{(a)_{\mu}}{(b)_{\mu}} = \prod_{i=1}^n \frac{(at^{1-i})_{\infty}}
{(bt^{1-i})_{\infty}} \sum_{\tau} \frac{(b/a)_{\tau}}{h'_{\tau}}
(at^{1-n})^{|\tau|} u_{\mu}\left(P_{\tau}(w)\right)
$$
It follows that
\begin{equation} \label{pana}
S = \prod_{i=1}^n \frac{(at^{1-i})_{\infty}}{(bt^{1-i})_{\infty}} 
\sum_{\tau,\sigma,\mu} \frac{(c)_{\sigma} z^{|\sigma|}}{h'_{\sigma}}
f^{\mu}_{\lambda\sigma}(q,t) x^{|\mu|} 
(at^{1-n})^{|\tau|} \frac{(b/a)_{\tau}}{h'_{\tau}} \: 
P_{\tau}(t^{\delta})\,u_{\tau}\left(P_{\mu}(w)\right) 
\end{equation}
where we have also used the symmetry property (\ref{symm.mac}) for Macdonald
polynomials. Now, from the definition of the coefficients
$f^{\mu}_{\lambda\sigma}(q,t)$ and the $q$-binomial theorem it also 
follows that
$$
\sum_{\sigma,\mu} \frac{(c)_{\sigma}z^{|\sigma|}}{h'_{\sigma}}
f^{\mu}_{\lambda\sigma}(q,t) x^{|\mu|}\, u_{\tau}\left(P_{\mu}(w) 
\right) = u_{\tau} \left( \prod_{i=1}^n\frac{(czxw_i)_{\infty}}
{(zxw_i)_{\infty}} x^{|\lambda|} P_{\lambda}(w) \right)
$$
Using this result in (\ref{pana}) gives 
\begin{eqnarray}
S &=& \prod_{i=1}^n \frac{(at^{1-i})_{\infty}(czxt^{n-i})_{\infty}}
{(bt^{1-i})_{\infty}(zxt^{n-i})_{\infty}} \sum_{\tau} 
\frac{(at^{1-n})^{|\tau|} (b/a)_{\tau} (zxt^{n-1})_{\tau}}
{h'_{\tau} (czxt^{n-1})_{\tau}} \,x^{|\lambda|} P_{\tau}(t^{\delta}) 
u_{\tau}\left(P_{\lambda}(w)\right) \nonumber\\
&=& \prod_{i=1}^n \frac{(at^{1-i})_{\infty}(czxt^{n-i})_{\infty}}
{(bt^{1-i})_{\infty}(zxt^{n-i})_{\infty}} x^{|\lambda|} 
P_{\lambda}(t^{\delta})
\sum_{\tau} \frac{(b/a)_{\tau} (zxt^{n-1})_{\tau}}{h'_{\tau}
(czxt^{n-1})_{\tau}} (at^{1-n})^{|\tau|} \,u_{\lambda}
\left(P_{\tau}(w)\right) \label{palta}
\end{eqnarray}
where again, we have used the symmetry property (\ref{symm.mac}). The
crucial point now is that if $b=aczxt^{n-1}$, then the sum over $\tau$  
can be carried out using the $q$-binomial theorem again, viz.
$$
u_{\lambda}\left( \frac{(zxt^{n-1})_{\tau}}{h'_{\tau}} P_{\tau}
(at^{1-n}w) \right) = \prod_{i=1}^n \frac{(azxt^{n-i})_{\infty}}
{(at^{1-i})_{\infty}} \; \frac{(a)_{\lambda}}{(azxt^{n-1})_{\lambda}}
$$
Using this in (\ref{palta}) gives the required result. \hfill $\Box$\\[2mm]
{\it Remarks}\\
{\it 1.}
When $\lambda=0$ this identity reduces to Gauss's formula (\ref{gauss.n}),
and can thus be considered a ``shifted'' version of Gauss's theorem.\\
{\it 2.}
The identity (\ref{chichi}) is simply the case $b=c=0$, $x=1$
of the above lemma.

\subsection{Sear's ${}_4\Phi_3$ transformation}

A final application of the $q$-Saalsch\"utz formulae (\ref{saal.n}) 
and (\ref{saal.n2}) is 
in deriving a multivariable analogue of Sear's transformation of a 
terminating, balanced ${}_4\phi_3$ series 
(see \cite[eq. (3.2.1)]{gasperbk}. Following the
proof given in \cite[Ex. 2.4]{gasperbk}, it suffices to write
the product of two specific ${}_2\Phi_1$ series in two different ways
and compare coefficients of the Macdonald polynomial $P_{(N^n)}(z)$ in 
each case. Indeed, the coefficient of $P_{(N^n)}(z)$ in the
product 
$$
P:={}_2\Phi_1\harg{a,b}{c}{z}\, {}_2\Phi_1\harg{d,e}{abde/c}{\frac{abz}{c}}
$$
can been seen to be
\begin{equation} \label{buscar.1}
\frac{(a)_{(N^n)} (b)_{(N^n)}}{(c)_{(N^n)}
h'_{(N^n)}}\; {}_4\Phi_3\harg{q^{-N},t^{n-1}q^{1-N}c^{-1},d,e}
{t^{n-1}q^{1-N}a^{-1},t^{n-1}q^{1-N}b^{-1},abde/c}{qt^{\delta}}
\end{equation}
where one must use the expression for $f^{(N^n)}_{\mu\nu}(q,t)$ given
in (\ref{imoto}), along with (\ref{ototo}) and (\ref{shumi}). On the
other hand we have by Euler's transformation (\ref{euler.2})
\begin{equation} \label{happy}
P = \prod_{i=1}^n \frac{(z_i)_{\infty}}{(abz_i/c)_{\infty}}\,
{}_2\Phi_1\harg{abe/c,abd/c}{abde/c}{z}\, {}_2\Phi_1\harg{a,b}{c}{z}
\end{equation}
Expanding this using the $q$-binomial theorem and then examining the
coefficient of $P_{(N^n)}(z)$, it follows from use of (\ref{imoto}),
(\ref{saal.n}), (\ref{ototo}) and (\ref{shumi}) that the coefficient 
of $P_{(N^n)}(z)$ in (\ref{happy}) is
\begin{equation} \label{buscar.2}
\left(\frac{ab}{c}\right)^{nN} \frac{(c/a)_{(N^n)} (c/b)_{(N^n)}}
{(c)_{(N^n)} h'_{(N^n)}} \; {}_4\Phi_3\harg{q^{-N},t^{n-1}q^{1-N}c^{-1},
abe/c,abd/c}{t^{n-1}q^{1-N}ac^{-1},t^{n-1}q^{1-N}bc^{-1},
abde/c}{qt^{\delta}}
\end{equation}
Comparing (\ref{buscar.1}), (\ref{buscar.2}), yields the multivariable
analogue of Sear's transformation
\begin{eqnarray}
{}_4\Phi_3\harg{q^{-N},t^{n-1}q^{1-N}c^{-1},d,e}
{t^{n-1}q^{1-N}a^{-1},t^{n-1}q^{1-N}b^{-1},abde/c}{qt^{\delta}} =
\hspace{6cm}\nonumber\\
\left(\frac{ab}{c}\right)^{nN} \frac{(c/a)_{(N^n)} (c/b)_{(N^n)}}
{(a)_{(N^n)} (b)_{(N^n)}} \; {}_4\Phi_3\harg{q^{-N},t^{n-1}q^{1-N}c^{-1},
abe/c,abd/c}{t^{n-1}q^{1-N}ac^{-1},t^{n-1}q^{1-N}bc^{-1},
abde/c}{qt^{\delta}}
\end{eqnarray}

\setcounter{equation}{0}
\section{Bilateral series}
We now turn our attention to bilateral basic hypergeometric series. 
The key point is that by using the following property of the
Macdonald polynomials \cite{mac}
\begin{equation} \label{shift}
|x|^a\,P_{\lambda}(x) = P_{\lambda+a}(x), \qquad a\in\zz, 
\qquad |x|:=x_1\cdots x_n
\end{equation}
where $\lambda+a := (\lambda_1+a,\ldots,\lambda_n+a)$, one can define 
Macdonald polynomials for all $n$-tuples $\lambda$ with
$\lambda_1\geq \lambda_2 \geq \cdots \geq \lambda_n$ with $\lambda_i\in\zz$.
namely if $\lambda_n<0$ then $P_{\lambda}(x) = |x|^{\lambda_n}\,
P_{\lambda-\lambda_n}(x)$. Denote the set of such partitions as ${\cal P}$
and those with non-negative entries as ${\cal P}_+$.

It also follows from (\ref{shift}) and the definition of the inner
product on the space of Macdonald polynomials that for any $\lambda\in
{\cal P}_+$ and any $a$ we have $\inn{|x|^a P_{\lambda}}{P_{\mu}}
=\inn{P_{\lambda}}{|x|^{-a} P_{\mu}}$ and it follows that we can extend
the definition of the inner product to polynomials $P_{\lambda}$ for
all $\lambda\in{\cal P}$, and that they remain orthogonal viz,
$$
\inn{P_{\lambda}}{P_{\mu}} = \delta_{\lambda\mu}\inn{P_{\lambda-\lambda_n}}
{P_{\lambda-\lambda_n}}
$$

Following Kaneko \cite{kaneko97a}, define the multivariable bilateral 
hypergoemetric series by
\begin{eqnarray}
{}_r\Psi_{s+1}\harg{a_1,\ldots,a_r}{b,b_1,\ldots,b_s}{z} :=
\hspace{8cm}\nonumber\\
\prod_{i=1}^n \frac{(bt^{i-1})_{\infty}(q)_{\infty}}
{(qt^{i-1})_{\infty}(b)_{\infty}}
\sum_{\lambda\in{\cal P}} \left((-1)^{|\lambda|} q^{n(\lambda')}
\right)^{s+1-r} \frac{(qt^{n-1})_{\lambda}(a_1)_{\lambda}\cdots 
(a_r)_{\lambda}}{(bt^{n-1})_{\lambda} (b_1)_{\lambda} \cdots (b_s)_{\lambda}
h'_{\lambda}} P_{\lambda}(z)
\end{eqnarray}
For $\lambda_n<0$, the above factorial symbols $(a)_{\lambda}$ must be
interpreted according to
$$
(a)_{\lambda} := t^{\sum_i (i-1)\lambda_i} \; \prod_{i=1}^n
\frac{(at^{1-i})_{\infty}}{(aq^{\lambda_i}t^{1-i})_{\infty}}
$$
and $h'_{\lambda}$ is interpreted according to (\ref{lapiz}).

\subsection{Kaneko's ${}_1\Psi_1$ summation formula}
In \cite{kaneko97a}, Kaneko gave a simple argument for the summation of
a ${}_1\Psi_1$ series. Here we give an alternative proof following
Andrew's argument in the one-variable case \cite{andr69a} which uses
Gauss's theorem. 

\begin{thm}[Kaneko]
For $|b/a|<|x_i|<1$ for all $i=1,\ldots,n$ we have
\begin{equation} \label{raman}
{}_1\Psi_1\harg{a}{b}{x} =\prod_{i=1}^n\frac{(ax_i)_{\infty}(q/ax_i)_{\infty}
(bt^{i-1}/a)_{\infty}(q)_{\infty}}
{(x_i)_{\infty}(b/ax_i)_{\infty} (qt^{i-1}/a)_{\infty}(b)_{\infty}}
\end{equation}
\end{thm}
{\bf Proof.} \quad We have from the $q$-binomial theorem that
$$
\prod_{i=1}^n \frac{(b/ax_i)_{\infty}}{(q/ax_i)_{\infty}} \;
{}_1\Psi_1\harg{a}{b}{x} = \prod_{i=1}^n \frac{(bt^{i-1})_{\infty}
(q)_{\infty}}{(qt^{i-1})_{\infty}(b)_{\infty}} 
\sum_{\stackrel{\sigma\in{\cal P}_{+}}{\mu\in{\cal P}}}
\frac{(b/q)_{\sigma}}{h'_{\sigma}}\left(\frac{q}{a}\right)^{|\sigma|}
\; \frac{(qt^{n-1})_{\mu}}{(bt^{n-1})_{\mu} h'_{\mu}}
P_{\sigma}(x^{-1}) P_{\mu}(x)
$$
We need to expand $P_{\sigma}(x^{-1}) P_{\mu}(x)$ in terms of
Macdonald polynomials. Suppose
$$
P_{\sigma}(x^{-1}) P_{\mu}(x) = \sum_{\lambda\in{\cal P}}
d^{\lambda}_{\sigma\mu} P_{\lambda}(x)
$$
{}From the orthogonality of $\{P_{\lambda}\}_{\lambda\in{\cal P}}$ we
have that 
\begin{eqnarray*}
\inn{P_{\nu-\nu_n}(x)}{P_{\nu-\nu_n}(x)} d^{\nu}_{\sigma\mu} &=& 
\inn{P_{\mu}(x)}{P_{\sigma}(x^{-1}) P_{\mu}(x)} \\
&=& \inn{P_{\nu-\nu_n}(x) P_{\sigma}(x)}{P_{\mu-\nu_n}(x)}
\end{eqnarray*}
Since both $\sigma$ and $\nu-\nu_n\in {\cal P}_+$, then $d^{\nu}_{\sigma
\mu}$ is non-zero only for all $\mu$ such that $\mu-\nu_n \in {\cal P}_+$.
Thus 
$$
d^{\nu}_{\sigma\mu} = \frac{\inn{P_{\mu-\nu_n}}{P_{\mu-\nu_n}}}
{\inn{P_{\nu-\nu_n}}{P_{\nu-\nu_n}}}\;f^{\mu-\nu_n}_{\nu-\nu_n,\sigma}
= \frac{h'_{\mu-\nu_n} h_{\nu-\nu_n}}{h_{\mu-\nu_n} h'_{\nu-\nu_n}}
\frac{(t^n)_{\mu-\nu_n} (qt^{n-1})_{\nu-\nu_n}}{(qt^{n-1})_{\mu-\nu_n}
(t^n)_{\nu-\nu_n}}\; f^{\mu-\nu_n}_{\nu-\nu_n,\sigma}
$$
Thus we have
\begin{equation} \label{calor}
\prod_{i=1}^n \frac{(b/ax_i)_{\infty}(b)_{\infty}}{(q/ax_i)_{\infty}
(q)_{\infty}} \; {}_1\Psi_1\harg{a}{b}{x} =
\prod_{i=1}^n \frac{(bt^{i-1})_{\infty}}{(qt^{i-1})_{\infty}}
\sum_{\nu\in{\cal P}} \frac{h_{\nu-\nu_n} (qt^{n-1})_{\nu-\nu_n}}
{h'_{\nu-\nu_n} (t^n)_{\nu-\nu_n}} \; B_{\nu} P_{\nu}(x)
\end{equation}
where 
$$
B_{\nu} = \prod_{i=1}^n \frac{(bq^{\nu_n}t^{n-i})_{-\nu_n}}
{(aq^{\nu_n}t^{1-i})_{-\nu_n}} \sum_{\sigma\in{\cal P}_+} \sum_{
{\mu\; {\rm s.t.} \atop \mu-\nu_n\in{\cal P}_+}}
\frac{(b/q)_{\sigma}}{h'_{\sigma}}
\left(\frac{q}{a}\right)^{|\sigma|} \frac{(aq^{\nu_n})_{\mu-\nu_n}}
{(bq^{\nu_n}t^{n-1})_{\mu-\nu_n}}\; P_{\mu-\nu_n}(t^{\delta})\,
f^{\mu-\nu_n}_{\nu-\nu_n,\sigma}(q,t)
$$
We now observe that we can use Lemma \ref{gauss.shift} to sum the above  
expression yielding 
$$
B_{\nu} = \prod_{i=1}^n \frac{(bq^{\nu_n}t^{n-i})_{-\nu_n}
(bt^{i-1}/a)_{\infty} (q^{\nu_n+1} t^{n-i})_{\infty}}
{(aq^{\nu_n}t^{1-i})_{-\nu_n}(qt^{i-1})_{\infty} (q t^{n-i}/a)_{\infty}}\; 
\frac{(aq^{\nu_n})_{\nu-\nu_n}}
{(q^{\nu_n+1}t^{n-1})_{\nu-\nu_n}}\; P_{\nu-\nu_n}(t^{\delta})
$$
Clearly  $B_{\nu}=0$ when $\nu_n<0$. In the cases when $\nu_n\geq 0$
(i.e. $\nu\in{\cal P}_+$) it follows after some manipulations that
$$
B_{\nu} = \prod_{i=1}^n \frac{(bt^{i-1}/a)_{\infty}(qt^{i-1})_{\infty}}
{(qt^{i-1}/a)_{\infty} (bt^{i-1})_{\infty}} 
\frac{h'_{\nu-\nu_n} (t^n)_{\nu-\nu_n}}{h_{\nu-\nu_n}
(qt^{n-1})_{\nu-\nu_n}} \frac{(a)_{\nu}}{h'_{\nu}}
$$
Substituting back into (\ref{calor}) and using the $q$-binomial
theorem gives the result. \hfill $\Box$

\subsection{${}_2\Psi_2$ transformations}

Our final application of the theory of Macdonald polynomials will be to 
transformation and summation formulae for bilateral ${}_2\Psi_2$
series. We begin with the following important result of Kadell and Kaneko
\cite{kaneko96a,kad97a,kaneko97b} which is equivalent to 
the ${}_1\Psi_1$ summation theorem given above

\begin{thm} \label{ct}
Given a partition $\lambda\in{\cal P}$ let
$$
A_{\lambda}(a,b) := {\rm C.T.}\left\{P_{\lambda}(x) \prod_{i=1}^n
(x_i)_a (q/x_i)_b  \Delta_q(x) \right\}
$$
where $\Delta_q(x)$ is the Macdonald weight function given in (\ref{dedo})
and $a$, $b$ are arbitrary complex numbers. Then 
$$
A_{\lambda}(a,b) = q^{(1+b)|\lambda|}\,\prod_{i=1}^n \frac{
(q^{1+a}t^{i-1})_{\infty} (q^{1+b}t^{i-1})_{\infty}}{
(qt^{i-1})_{\infty} (q^{1+a+b}t^{i-1})_{\infty}}\; \frac{
(q^{-b})_{\lambda}}{(q^{1+a}t^{n-1})_{\lambda}}\;P_{\lambda}(t^{\delta})
\inn{1}{1}
$$
\end{thm}
The connection with the ${}_1\Psi_1$ formula is as follows: suppose
\begin{equation} \label{zapato.1}
\prod_{i=1}^n (x_i)_a (q/x_i)_b = \sum_{\lambda\in{\cal P}} c_{\lambda}
(a,b)\,P_{\lambda}(x^{-1})
\end{equation}
It follows from the orthogonality of Macdonald polynomials that 
$c_{\lambda}(a,b) = A_{\lambda}(a,b)/\inn{P_{\lambda}}{P_{\lambda}}$
with $A_{\lambda}(a,b)$ given as in Theorem \ref{ct}. Now, the
transformation $x_i\rightarrow q^{-u}/x_i$ certainly doesn't affect the
constant term of the expression appearing in Theorem \ref{ct}. Thus
\begin{equation} \label{llave}
{\rm C.T.}\left\{P_{\lambda}(x^{-1}) \prod_{i=1}^n
(q^{-u}/x_i)_{a'} (q^{1+u}x_i)_{b'}  \Delta_q(x) \right\}
= q^{u|\lambda|}\, A_{\lambda}(a',b')
\end{equation}
As above, it follows that if
\begin{equation} \label{zapato.2}
\prod_{i=1}^n (q^{-u}/x_i)_{a'} (q^{1+u}x_i)_{b'} = 
\sum_{\mu\in{\cal P}} c_{\mu}' (a',b')\,P_{\mu}(x)
\end{equation}
then $c_{\mu}' (a',b')=q^{u|\mu|} A_{\mu}(a',b')/\inn{P_{\mu}}{P_{\mu}}$.
Define
\begin{equation} \label{idef}
I(a,b,a',b';u) := {\rm C.T.} \left\{ \prod_{i=1}^n (x_i)_a (q/x_i)_b
(q^{-u}/x_i)_{a'} (q^{u+1} x_i)_{b'} \,\Delta_q(x) \right\}
\end{equation}
It follows from (\ref{zapato.1}) and (\ref{zapato.2}) that
\begin{eqnarray}
I(a,b,a',b';u) &=& \sum_{\lambda\in{\cal P}} c_{\lambda}(a,b)
c'_{\lambda}(a',b') \inn{P_{\lambda}}{P_{\lambda}} \nonumber\\
&=&\prod_{i=1}^n \frac{(q^{1+a})_{\infty} (q^{1+b}t^{i-1})_{\infty}
(q^{1+a'}t^{i-1})_{\infty} (q^{1+b'}t^{i-1})_{\infty}} 
{(q)_{\infty} (qt^{i-1})_{\infty}
(q^{1+a+b}t^{i-1})_{\infty} (q^{1+a'+b'}t^{i-1})_{\infty}} \inn{1}{1}
\times\nonumber\\
&&\times {}_2\Psi_2\harg{q^{-b},q^{-b'}}{q^{1+a},q^{1+a'}t^{n-1}}
{q^{b+b'+2+u} t^{\delta}} \label{manzana}
\end{eqnarray}
We now note there are a couple of values of $u$ for which we can
derive an alternative expression for $I(a,b,a',b';u)$.
\bigskip\\
{\bf The case $u=a-1$}\\
In this case, note that $(x_i)_a (q^a x_i)_{b'} = (x_i)_{a+b'}$.
Using the $q$-binomial theorem to expand $\prod_i(q^{-u}/x_i)_{a'}$
and inserting in (\ref{idef}) we have
\begin{eqnarray}
I &=& \sum_{\lambda\in{\cal P}_+} \frac{(q^{-a'})_{\lambda}}{h'_{\lambda}}
q^{(a'+1-a)|\lambda|} {\rm C.T.}\left\{ \prod_{i=1}^n (x_i)_{a+b'}
(q/x_i)_b P_{\lambda}(x^{-1}) \Delta_q(x) \right\} \nonumber\\
&=& \prod_{i=1}^n \frac{(q^{1+b}t^{i-1})_{\infty} (q^{1+a+b'}t^{i-1})_{
\infty}}{(qt^{i-1})_{\infty} (q^{1+a+b+b'}t^{i-1})_{\infty}} \inn{1}{1}
{}_2\Phi_1 \harg{q^{-a'},q^{-a-b'}}{q^{1+b}t^{n-1}}{q^{1+a'+b'} t^{\delta}}
\label{hijastra}
\end{eqnarray}
where we have used (\ref{llave}) with $u=-1$. Comparing (\ref{hijastra})
with the $u=a-1$ case of (\ref{manzana}) gives (upon setting
$\alpha_1= q^{-b}$, $\alpha_2= q^{-b'}$, $\alpha_3= q^{1+a}$ and
$\alpha_4= q^{1+a'}$) the transformation formula
\begin{eqnarray}
{}_2\Psi_2\harg{\alpha_1,\alpha_2}{\alpha_3,\alpha_4 t^{n-1}}
{\frac{\alpha_3}{\alpha_1\alpha_2}\,t^{\delta}} =
\prod_{i=1}^n \frac{(q)_{\infty} (\alpha_3 t^{i-1}/\alpha_1)_{\infty}
(\alpha_4 t^{i-1}/\alpha_2)_{\infty}(\alpha_3 t^{i-1}/\alpha_2)_{\infty}}
{(\alpha_3)_{\infty} (\alpha_4 t^{i-1})_{\infty}(qt^{i-1}/\alpha_2)_{\infty}
(\alpha_3 t^{i-1}/\alpha_1\alpha_2)_{\infty}} \times\nonumber\\
\times\; {}_2\Phi_1\harg{q/\alpha_4, q\alpha_2/\alpha_3}{qt^{n-1}/\alpha_1}
{\frac{\alpha_4}{\alpha_2}\,t^{\delta}} \hspace{3cm} \label{beetroot}
\end{eqnarray}
Note that when $\alpha_3 = q\alpha_1$, the ${}_2\Phi_1$ appearing in
the above equation can be summed by Gauss's formula (\ref{gauss.n}) 
yielding the summation formula
\begin{equation}
{}_2\Psi_2\harg{\alpha_1,\alpha_2}{q\alpha_1,\alpha_4 t^{n-1}}
{\frac{q}{\alpha_2} t^{\delta}} = \prod_{i=1}^n \frac{(q)_{\infty}
(qt^{i-1})_{\infty} (q\alpha_1 t^{i-1}/\alpha_2)_{\infty}
(\alpha_4 t^{i-1}/\alpha_1)_{\infty}}
{(q\alpha_1)_{\infty} (\alpha_4 t^{i-1})_{\infty} (qt^{i-1}/\alpha_2)_{
\infty} (qt^{i-1}/\alpha_1)_{\infty}}
\end{equation}
\bigskip
{\bf The case $u=-1$}\\
Note from the definition (\ref{idef}) that
$$
I(a,b,a',b';-1) = I(b',b,a',a;-1) = I(a,a',b,b';-1)
$$
Taking for example, the first of these relations, and using the
equivalent expressions given by (\ref{manzana}) gives the following
transformation formula
\begin{eqnarray}
{}_2\Psi_2\harg{\alpha_1,\alpha_2}{\alpha_3,\alpha_4 t^{n-1}}
{\frac{q}{\alpha_1\alpha_2} t^{\delta}} = \prod_{i=1}^n
\frac{(q/\alpha_2)_{\infty} (\alpha_3 t^{i-1})_{\infty}
(\alpha_4 t^{i-1}/\alpha_2)_{\infty} (\alpha_3 t^{i-1}/\alpha_1)_{\infty}}
{(\alpha_3)_{\infty} (q t^{i-1}/\alpha_2)_{\infty}
(q t^{i-1}/\alpha_1\alpha_2)_{\infty} (\alpha_3\alpha_4 t^{i-1}/q)_{\infty}}
\times \nonumber\\
\times\; {}_2\Psi_2\harg{\alpha_1,q\alpha_3}{q/\alpha_2,\alpha_4 t^{n-1}}
{\frac{\alpha_3}{\alpha_1} t^{\delta}} \hspace{4cm} \label{corn}
\end{eqnarray}

\subsection{The $q\rightarrow 1$ limit}

One of the curious features of the bilateral basic series ${}_2\psi_2$
is that there appears to be no analogue of Gauss's summation for
the ${}_2\phi_1$ series, whereas in the $q=1$ case, the bilateral
${}_2H_2$ series with unit argument can indeed be summed \cite{doug07}.
Let us now show that this behaviour carries over in the multivariable
case. We begin by defining the $q\rightarrow 1$ limit of the bilateral
series ${}_r\Psi_{s+1}$, namely
\begin{equation}
{}_r H_{s+1}\harg{a_1\ldots,a_r}{b,b_1\ldots,b_s}{z}
:= \prod_{i=1}^n \frac{\Gamma(1+k(i-1)) \Gamma(b)}
{\Gamma(b+k(i-1))} \sum_{\lambda\in{\cal P}} \frac{(1+k(n-1))_{\lambda}
(a_1)_{\lambda}\cdots (a_r)_{\lambda}}{(b+k(n-1))_{\lambda}
(b_1)_{\lambda}\cdots (b_s)_{\lambda}d'_{\lambda}} \; P_{\lambda}(x;k)
\end{equation}
Here $P_{\lambda}(x;k)$ is the Jack polynomial (normalized so that the
coefficient of the monomial symmetric function $m_{\lambda}(x)$ is unity),
the generalized Pochhammer symbols are defined by
$(a)_{\lambda} := \prod_i (a+k(1-i))_{\lambda_i}$ with $(a)_n:=
a(a+1)\cdots (a+n-1)$ being the usual Pochhammer symbol, and
$d'_{\lambda} := \prod_{s\in\lambda}(a(a)+1+k.l(s))$ which is the
$q\rightarrow 1$ limit of $h'_{\lambda}$ defined in (\ref{libro}).

It follows from (\ref{manzana}) that 
\begin{eqnarray*}
\lim_{q\rightarrow 1} I(a,b,a',b';u) = \prod_{i=1}^n \frac{
\Gamma(1+k(i\!-\!1)) \Gamma(1+a+b+k(i\!-\!1)) \Gamma(1+a'+b'+k(i\!-\!1))}
{\Gamma(1+a) \Gamma(1+b+k(i\!-\!1)) \Gamma(1+b'+k(i\!-\!1)) 
\Gamma(1+a'+k(i\!-\!1))} \times \\
\times\; \inn{1}{1}_k\;{}_2 H_2 \harg{-b,-b'}{1+a,1+a'+k(n-1)}{1^n}
\end{eqnarray*}
where $\inn{1}{1}_k$ is the appropriate limit of the corresponding 
Macdonald quantity, given in (\ref{ground}). The crucial point now is that 
in the $q=1$ case,
we can combine the quantities $a$ and $b'$ (respectively $a'$ and $b$)
in the definition of $I(a,b,a',b';u)$ in (\ref{idef}), so that
\begin{eqnarray*}
\lim_{q\rightarrow 1} I(a,b,a',b';u) &=& \lim_{q\rightarrow 1}
A_{0}(a+b',a'+b) \\
&=& \prod_{i=1}^n \frac{ \Gamma(1+k(i-1)) \Gamma(1+a+a'+b+b'+k(i-1))}
{\Gamma(1+a+b'+k(i-1)) \Gamma(1+a'+b+k(i-1))} \inn{1}{1}_k
\end{eqnarray*}
{}From these two equations, we have the summation formula
\begin{eqnarray}
{}_2 H_2\harg{-b,-b'}{1+a,1+a'+k(n-1)}{1^n} = \prod_{i=1}^n\frac{\Gamma(1+a)
\Gamma(1+a+a'+b+b'+k(i-1))}{\Gamma((1+a+b+k(i-1)) \Gamma(1+a'+b'+k(i-1))}
\times\nonumber\\[3mm]
\times\; \prod_{i=1}^n\frac{ \Gamma(1+b+k(i-1)) \Gamma(1+b'+k(i-1)) 
\Gamma(1+a'+k(i-1))}{\Gamma(1+a+b'+k(i-1)) \Gamma(1+a'+b+k(i-1))}
\qquad
\end{eqnarray}
which certainly reduces to Dougall's result \cite{doug07} when $n=1$.

\subsection{Bailey's ${}_2\Psi_2$ transformations}

In this final section, we shall derive multivariable analogues of
Bailey's general transformations for ${}_2\psi_2$ series \cite{bail50a}.
First, given a general partition $\lambda\in {\cal P}$, define
$-\lambda^R := (-\lambda_n,-\lambda_{n-1},\ldots,-\lambda_1)$. It
follows from (\ref{sombrero}) (which can easily be seen to hold
for all partitions $\sigma\in{\cal P}$) and (\ref{shift}) that
\begin{equation}
P_{-\lambda^R}(z) = P_{\lambda}(z^{-1}), \qquad \forall \lambda\in{\cal P}
\end{equation}
It also follows from the definitions that
\begin{equation} \label{jb}
(a)_{-\lambda^R} = \frac{(-q/a)^{|\lambda|} q^{n(\lambda')}}
{(qt^{n-1}/a)_{\lambda}} \hspace{3cm}
\frac{(qt^{n-1})_{-\lambda^R}}{h'_{-\lambda^R}} = t^{(1-n)|\lambda|}
\frac{(qt^{n-1})_{\lambda}}{h'_{\lambda}}
\end{equation}

The proof follows \cite{bail50a}: from the ${}_1\Psi_1$ summation formula
(\ref{raman}) it follows that
\begin{eqnarray}
{}_1\Psi_1\harg{b}{c}{az}\; {}_1\Psi_1\harg{ab'}{c'}{z} =
\prod_{i=1}^n \frac{(ct^{i-1}/b)_{\infty} (c't^{i-1}/ab')_{\infty}
(qt^{i-1}/b')_{\infty} (qt^{i-1}/ab)_{\infty}}
{(qt^{i-1}/b)_{\infty} (qt^{i-1}/ab')_{\infty} (c't^{i-1}/b')_{\infty}
(ct^{i-1}/ab)_{\infty}}  \times \nonumber\\
\times\; {}_1\Psi_1\harg{b'}{c'}{az} {}_1\Psi_1\harg{ab}{c}{z}
\label{nectarine}
\end{eqnarray}
We now need to compare the constant terms on each sides. Certainly, the
product $P_{\lambda}(z)\,P_{\mu}(z)$ only has a constant term if
$\mu=-\lambda^R$ and that constant term is $\inn{P_{\lambda}}{P_{\lambda}}/
\inn{1}{1}$. If we now use this fact to write down the constant terms
in each side of (\ref{nectarine}), using the simplifications (\ref{jb})
and changing parameters, we end up with the transformation
\begin{eqnarray*}
{}_2\Psi_2 \harg{\alpha_1,\alpha_2}{\alpha_3,\alpha_4 t^{n-1}}{z t^{\delta}}
= \prod_{i=1}^n \frac{ (\alpha_3 t^{i-1}/\alpha_1)_{\infty}
(\alpha_4 t^{i-1}/\alpha_2)_{\infty} (\alpha_2 z t^{i-1})_{\infty}
(q\alpha_4 t^{i-1}/\alpha_1\alpha_2 z)_{\infty}}
{(qt^{i-1}/\alpha_1)_{\infty} (\alpha_4 t^{i-1})_{\infty}
(zt^{i-1})_{\infty} (\alpha_3\alpha_4 t^{i-1}/\alpha_1\alpha_2 z)_{\infty}}
\times \\
\times\;\prod_{i=1}^n \frac{(\alpha_3 t^{i-1})_{\infty} (q/\alpha_2)_{\infty}}
{(qt^{i-1}/\alpha_2)_{\infty} (\alpha_3)_{\infty}}\;
{}_2\Psi_2 \harg{q/\alpha_2 z, q/\alpha_3}{q/\alpha_2,qt^{n-1}\alpha_4/
\alpha_1\alpha_2 z}{\frac{\alpha_3}{\alpha_1} t^{\delta}}
\end{eqnarray*}
Following \cite{bail50a}, we can rewrite this using the identity
$$
{}_2\Psi_2\harg{a,b}{c,d}{z} = \prod_{i=1}^n \frac{(ct^{i-1})_{\infty}
(q/a)_{\infty}}{(c)_{\infty} (qt^{i-1}/a)_{\infty}}\;
{}_2\Psi_2\harg{q/c,qt^{n-1}/d}{q/a,qt^{n-1}/b}{\frac{cd}{ac}\:z^{-1}}
$$
which follows from making the substitution $\lambda\rightarrow -\lambda^R$
in the definition of the ${}_2\Psi_2$ series and simplifying using 
(\ref{jb}). The final result is thus
\begin{eqnarray}
{}_2\Psi_2\harg{\alpha_1,\alpha_2}{\alpha_3,\alpha_4 t^{n-1}}{z t^{\delta}}
= \prod_{i=1}^n \frac{(\alpha_3 t^{i-1})_{\infty} (\alpha_3 t^{i-1}/\alpha_1
)_{\infty} (\alpha_4 t^{i-1}/\alpha_2)_{\infty} (\alpha_2 z)_{\infty}
(q\alpha_4t^{i-1}/\alpha_1\alpha_2 z)_{\infty}}
{(\alpha_3)_{\infty} (qt^{i-1}/\alpha_1)_{\infty} (\alpha_4 t^{i-1})_{\infty}
(z t^{i-1})_{\infty} (\alpha_3\alpha_4 t^{i-1}/\alpha_1\alpha_2 z)_{\infty}}
\times \nonumber\\
\times\; {}_2\Psi_2\harg{\alpha_2,\alpha_1\alpha_2 z/\alpha_4}
{\alpha_2 z, \alpha_3 t^{n-1}}{\frac{\alpha_4 t^{n-1}}{\alpha_2} t^{\delta}}
\hspace{3cm} \label{bt}
\end{eqnarray}
We remark here that the transformations appearing in previous section
(\ref{beetroot}), (\ref{corn}) are special cases of the transformation
(\ref{bt}).

\vspace{5mm}\noindent
{\bf\Large Acknowledgements}\\[2mm]
The authors thank George Andrews for discussions. They acknowledge 
the financial support of the Australian Research Council.

\bibliographystyle{plain}
%\bibliography{symm}

\end{document}